# Optimal control of differential inclusions with endpoint and state constraints and duality

**Elimhan N. Mahmudov**
Department of Mathematics, Istanbul Technical University, Istanbul, Turkey,
Azerbaijan National Academy of Sciences Institute of Control Systems, Baku, Azerbaijan.
elimhan22@yahoo.com, ORCID ID: 0000-0003-2879-6154

**Abstract** The paper studies optimal control problem described by higher order evolution differential inclusions (DFIs) with endpoint and state constraints. In the term of Euler-Lagrange type inclusion is derived sufficient condition of optimality for higher order DFIs. It is shown that the adjoint inclusion for the first order DFIs, defined in terms of locally adjoint mapping, coincides with the classical Euler-Lagrange inclusion. Then a duality theorem is proved, which shows that Euler-Lagrange inclusions are "duality relations" for both problems. At the end of the paper duality problems for third order linear and fourth order polyhedral DFIs are considered.

## 1. Introduction

The paper deals with the Mayer problem of the higher-order evolution differential inclusions (DFIs) with both endpoint and state constraints:

$$\text{minimize } f\left(x(0), x(T)\right), \tag{1}$$

(PHC)
$$\frac{d^k x(t)}{dt^k} \in F(x(t),t), \text{ a.e. } t \in [0,T], \tag{2}$$

$$\left(x^{(j)}(0), x^{(j)}(T)\right) \in S,\ j = 0,1,\ldots,k-1, \tag{3}$$

$$x(t) \in X(t),\ \forall t \in [0,T], \tag{4}$$

where $F(\cdot,t): \mathbb{R}^n \rightrightarrows \mathbb{R}^n$ is a convex set-valued mapping, $f: \mathbb{R}^{2n} \to \mathbb{R}^1$ is a proper convex function, $k$ is an arbitrary fixed natural number, $S \subseteq \mathbb{R}^{2n}$ and $X(t) \subseteq \mathbb{R}^n,\ \forall t \in [0,T]$ are nonempty subsets, $T$ is an arbitrary positive real number. The problem is to find an arc $\tilde{x}(\cdot)$ of the problem (PHC) that almost everywhere (a.e.) satisfies the inclusion (2) on the time interval $[0,T]$, everywhere constraints (3), (4), and minimizes the so-called Mayer's functional $f\left(x(0), x(T)\right)$. Let us refine the definition of the concept of a solution of problem for $k$-th order DFIs (2)-(4); suppose $AC^j\left([0,T], \mathbb{R}^n\right)$ is the space of $j$-times differentiable functions $x(\cdot): [0,T] \to \mathbb{R}^n$, where $j$-th order derivative $x^{(j)}(\cdot) \equiv d^{(j)}x(\cdot)/dt^j$ $(j = 1,\ldots,k)$ is absolutely continuous and $L^1\left([0,T], \mathbb{R}^n\right)$ is the Banach space of integrable (in the Lebesgue sense) functions $u(\cdot): [0,T] \to \mathbb{R}^n$ endowed with the norm $\|u(\cdot)\|_1 = \int_0^T |u(t)|dt.$ A function $x(\cdot) \in AC^{k-1}\left([0,T], \mathbb{R}^n\right)$ is a feasible solution of a problem (2)-(4) if there exists an integrable function $\upsilon(t) \in F(x(t),t),$ a.e. $t \in [0,T]$, with $x^{(k)}(t) = \upsilon(t)$ a.e. $t \in [0,T]$ and $x(\cdot)$ satisfies conditions (3),(4).
It should be noted that having endpoint constraint $(x(0), x(T)) \in S$ and state constraint $x(t) \in X(t), \ \forall t \in [0, T]$ as aspects of the model are required for applications. For example, Clarke [3, Chap. III] gives an excellent introduction to this problem with first order DFIs and describes several applications.

The present problem (PHC) in a nonconvex setting is an essential generalization of the problem considered in the paper of Loewen and Rockafellar [15] in which under several stringent conditions the necessary conditions for optimality are derived; such restrictions are imposed as constraint qualification; locally Lipschitz property of the cost functional $f$; measurability of $F$; the closedness of the values of $X(t)$ and the lower semicontinuity of $X(t)$, etc. Given the existing assumptions, a lot of effort has been made to overcome the significant difficulty that has arisen in the formulation and derivation of the necessary conditions. We believe that our optimality conditions contain more convenient forms of the transversality condition and Euler-Lagrange inclusion conditions, and therefore they are more practical than those previously published in [15] and [6]. Moreover, the simplicity of the LAM approach and the method of the “cone of tangent directions” instead of the normal cone simplifies the derivation and formulation of optimality conditions [18]. We hope that all these improvements will serve for the further development of the theory of duality theory.

We recall that in the theory of high-order differential inclusions there are many interesting papers on qualitative research, mainly related to existence theorems [1,2,4,5,7,8,10,11,13,14, 16,17,20-27,29-33]. In the paper [4] it is shown that for the optimal solution of the Mayer problem with a first-order differential inclusion there is a "dual variable" for which, in terms of the Hamiltonian function, a certain "dual relation" and transversality conditions are satisfied.
In the work of Frankowska [10] the considered classical optimal control problem of Bolza type with phase constraints is solved using the Hamilton-Jacobi equation. In the paper of Quincampoix [32] it examines where and how solutions related to differential inclusion may or may not fit into a given target. For this, sections of the target boundary are linked to the dynamics of the system. The behaviour of these solutions is qualitatively described in terms of the kernels of viability and the invariance of the sets. In the work [11] the authors investigate the viability problem for first-order differential inclusions. The paper [33] develops a tube theory for solving differential inclusions within a given set of sets. The concept is defined as the minimum invariant tube with values in the collection. Under certain requirements for the set, the existence and Lipschitz property of tubes with a solution is proved. In the paper [1] formulates the conditions under which there exists a solution to a viable problem with a second-order differential inclusion. In the paper [29], by extending the method of guiding functions, are given sufficient conditions for the existence of solutions to the second order where a set-valued map is not necessarily convex-valued. The paper [27] deals with the Mayer problem for first order differential inclusions and terminal functional constraints.

Although a significant part of the paper is devoted to the derivation of optimality conditions for problem (PHC) its main goal is to construct and study the duality theory for it. On the one hand, duality theory provides a powerful theoretical tool for the analysis of optimization and variational problems, and on the other hand, it opens the way to the development of new algorithms to solve them. The reader can refer to [7-9,12,25,28,34.35] and their references for more details on this topic. For convex optimization problems, the duality gap, that is, the difference between the optimal values of the primary and dual problems, is zero under a constraint qualification condition. As far as we know, there only are several papers [7-9,12,25,34,35] devoted to the problems of duality of ordinary and partial DFIs. In the work [9] using Fenchel-Lagrange duality is formulated conditions ensuring the strong duality

In the works [16-18,20-25], for optimal control problems with ordinary and partial DFIs in terms of locally adjoint mappings (LAMs) the optimality conditions are derived. The work [25] is devoted to optimization of second-order ordinary DFIs; the problems of viability of the second order for differential inclusions with end endpoint constraint are considered and duality theorems are proved.

The obtained results can be organized in the following order:

In *Section 2*, for the convenience of the readers, all definitions, basic facts and concepts from the book of Makhmudov [18] are given.

In *Section 3*, sufficient condition of optimality for a problem (PHC) with *k*-th order differential inclusion is proved. Also are formulated the so-called transversality conditions imposing some conditions on the endpoints of the trajectory $x^{*(j)}(0), x^{*(j)}(T),\ j=0,...,k-1.$ It is shown that in the particular case when $k=1$ the existing optimality conditions imply the classical Euler-Lagrange adjoint inclusion. Thus, we obtain sufficient optimality conditions for the problem posed by Loewen and Rockafellar [15]. At the end of the section, the results obtained in a linear third-order optimal control problem (PTL) are demonstrated. An analogue of the adjoint equation, transversality condition and Ponrtyagin's maximum principle are formulated.

*In Section 4* are investigated duality results for primary problem (PHC); for the optimality of the family of “dual variables” in the dual problem (PHC*) it is necessary and sufficient that the optimality conditions be satisfied. It simply means that the Euler-Lagrange inclusion is a “dual relation” for both the primary (PHC) and dual (PHC *) problem.

*In Section 5* duality problems for third order linear and fourth order polyhedral DFIs are investigated; first the dual problem is constructed for the linear problem of the third order, considered in Section 3. It is shown that for the function of the dual variable in the dual problem (PTL*) to be optimized, it is necessary and sufficient that the adjoint third order equation and the Pontryagin maximum principle be satisfied. Second, a dual problem is constructed for a fourth-order polyhedral problem. For this, the given problem is reduced to a linear programming problem, and thus the "support function" is calculated for the graph of the polyhedral mapping. It turns out that, according to convex programming problems, maximization in the dual problem is carried out over nonnegative functions. From an applied point of view, these examples show that the considered approach to constructing duality turns out to be justified.

In conclusion, we note that for the values of the primary ($\alpha$) and dual problems ($\alpha^*$), the inequality $\alpha \ge \alpha^*$ is satisfied, more precisely, if the variables of the primary and dual problems satisfy the Euler-Lagrange inclusion, then these values are equal. Moreover, it is clear that convex and convexified problems have the same convex dual problem. This is explained by the fact that convex, convex closed sets and convex hulls of nonconvex sets have the same support function. Here we study the optimality conditions for Mayer problems (PHC) and duality theorem based on dual operations of addition and infimal convolution of convex functions. However, the constructions of Euler-Lagrange type inclusions, transversality conditions, and duality problems are beyond the scope of this paper, so it is omitted. And in this sense, the results obtained in Sections 3 and 4 are only the visible part of the "icebergs".

## 2. Necessary facts, preliminaries

All definitions and concepts that we come across can be found in Mahmudov's book [18].

Suppose that $F:\mathbb{R}^n \rightrightarrows \mathbb{R}^n$ is a set-valued mapping from $n$-dimensional Euclidean space $\mathbb{R}^n$ into the family of subsets of $\mathbb{R}^n$, $\langle x,v\rangle$ be an inner product of $x$ and $v$. $F$ is convex closed if its graph is a convex closed set in $\mathbb{R}^{2n}$. Let's give important definitions, which we will often see in the paper:

$$H_F(x,v^*)=\sup_v\left\{\left\langle v,v^*\right\rangle : v\in F(x)\right\},\ v^*\in\mathbb{R}^n,$$

$$F_A(x;v^*)=\left\{v\in F(x):\left\langle v,v^*\right\rangle = H_F(x,v^*)\right\}.$$

$H_F$ and $F_A$ are called Hamiltonian function and argmaximum set for a set-valued mapping $F$, respectively. If $F(x)=\varnothing$ for a convex $F$ we put $H_F(x,v^*)=-\infty$.

Throughout this paper the support function of a set $Q\subseteq\mathbb{R}^n$ is denoted by

$$W_Q(x^*)=\sup\left\{\left\langle x,x^*\right\rangle : x\in Q\right\}.$$

For such a mapping $F$, the cone of tangent directions at the point $(x^0,v^0)\in\operatorname{gph}F$ is defined as follows

$$K_F(x^0,v^0)\equiv K_{\operatorname{gph}F}(x^0,v^0)=\operatorname{cone}\left[\operatorname{gph}F-(x^0,v^0)\right]=\left\{(\bar{x},\bar{v}):\right.$$
$$\left.\bar{x}=\gamma(x-x^0),\ \bar{v}=\gamma(v-v^0),\ \gamma>0\right\},\forall\ (x,v)\in\operatorname{gph}F.$$

A set-valued mapping $F^*(\cdot,x,v):\mathbb{R}^n \rightrightarrows \mathbb{R}^n$ defined by

$$F^*\left(v^*;(x^0v^0)\right)=\left\{x^*:(x^*,-v^*)\in K_F{}^*(x^0,v^0)\right\}$$

is called the LAM to $F$ at a point $(x^0,v^0)\in\operatorname{gph}F$, where $K_F{}^*(x^0,v^0)$ is the dual cone. Note that, using the definition of the cone of tangent vectors in the non-convex case, the LAM for non-convex multivalued mappings is determined by the same formula [18,p.129].

A "dual" mapping defined by

$$F^*(v^*;(x^0,v^0)):=\left\{x^*:H_F(x,v^*)-H_F(x^0,v^*)\right.$$
$$\left.\le\left\langle x^*,x-x^0\right\rangle,\ \forall x\in\mathbb{R}^n\right\},\ v\in F(x;v^*)$$

is called the LAM to "nonconvex" mapping $F$ at a point $(x^0,v^0)\in\operatorname{gph}F$. Obviously, in the convex case $H_F(\cdot,v^*)$ is concave and the latter definition of LAM coincide with the previous definition of LAM. In fact, the given in the paper notion LAM is closely related to the coderivative (Frechet) concept of Mordukhovich [27,28], $D^*F(x,v)(v^*)=\left\{x^*:(x^*,-v^*)\right.$ $\left.\in N\left((x,v);\operatorname{gph}F\right)\right\}$, which is essentially different for nonconvex mappings. Here $N((x,v);$ $\operatorname{gph}F)$ is a normal cone at $(x,v)\in\operatorname{gph}F$. In the most interesting settings for the theory and applications, coderivatives are nonconvex-valued and hence are not tangentially/derivatively generated. However, for the convex maps the two notions are equivalent.

A function $g=g(x,y)$ is called a proper function if it does not assume the value $-\infty$ and is not identically equal to $+\infty$. Obviously, $g$ is proper if and only if $\operatorname{dom}g\neq\varnothing$ and $g(x,y)$ is finite for $(x,y)\in\operatorname{dom}g=\left\{(x,y):g(x,y)<+\infty\right\}$.

**Definition 2.1** A function $g(x,y)$ is a closure if $\operatorname{epi}g=\left\{(\xi,x,y):\ \xi\ge g(x,y)\right\}$ is a closed set.

**Definition 2.2** The function $g^*\left(x^*,y^*\right)$ defined as below is called the conjugate of $g$:

$$g^*\left(x^*, y^*\right) = \sup_{x,y}\left\{\left\langle x, x^*\right\rangle + \left\langle y, y^*\right\rangle - g(x,y)\right\}.$$

It is easy to see that the function

$$M_G(x^*, v^*) = \inf\left\{\left\langle x, x^*\right\rangle - \left\langle v, v^*\right\rangle : (x,v) \in \operatorname{gph}F\right\} = \inf_x\left\{\left\langle x, x^*\right\rangle - H_F(x, v^*)\right\}.$$

is a support function taken with a minus sign. Besides, it follows that for a fixed $v^*$

$$M_F(x^*, v^*) = -\left[-H_F(\cdot, v^*)\right]^*(x^*)$$

that is, $M_F$ is the conjugate function for $-H_F(\cdot, v^*)$ taken with a minus sign. By Lemma 2.6 [18, p.64] it is noteworthy to see that $x^*$ is an element of the LAM $F^*$ if and only if

$$M_F(x^*, v^*) = \left\langle x, x^*\right\rangle - H_F(x, v^*).$$

## 3. Sufficient condition of optimality for a problem (PHC) with *k*-th order differential inclusion

In this section, we formulate a Euler-Lagrange inclusion for the problem under consideration. Due to the fact that the construction of the Euler-Lagrange inclusion, as well as transversality conditions are complicated by the accompaniment of discrete and discrete-approximation problems (see, for example [20-24]) we omit it and formulate only the final result. So let us formulate for a Lagrange problem (PHC) with $k$-th order differential inclusion the following Euler-Lagrange type inclusion and the so-called transversality conditions:

$$(a)\quad (-1)^k x^{*(k)}(t) \in F^*\left(x^*(t); (\tilde{x}(t), \tilde{x}^{(k)}(t)), t\right) + K^*_{X(t)}(\tilde{x}(t)),, \text{ a.e. } t \in [0,T]$$

$$(b)\quad \tilde{x}^{(k)}(t) \in F_A\left(\tilde{x}(t); x^*(t), t\right),$$

$$(c)\quad \left((-1)^{k-1} x^{*(k-1)}(0), (-1)^k x^{*(k-1)}(T)\right) \in \partial f(\tilde{x}(0), \tilde{x}(T)) - K_S^*(\tilde{x}(0), \tilde{x}(T)),$$

$$(d)\quad \left((-1)^{j+1} x^{*(j)}(0), (-1)^j x^{*(j)}(T)\right) \in K_S^*\left(\tilde{x}^{(k-j-1)}(0), \tilde{x}^{(k-j-1)}(T)\right),\ j = 0, \ldots, k-2.$$

The definition of the solution to the Euler - Lagrange inclusion is defined similarly to the definition of the solution to problem (PHC). A function $x^*(\cdot) \in AC^{k-1}\left([0,T],\ \mathbb{R}^n\right)$ is called a feasible solution of problem (*a*)-(*d*) if there exists a function $w(\cdot) \in L^1\left([0,T],\ \mathbb{R}^n\right)$ with $w(t) \in F^*\left(x^*(t); (x(t), x^{(k)}(t)), t\right) + K^*_{X(t)}(\tilde{x}(t))$, a.e. $t \in [0,T]$ such that $(-1)^k x^{*(k)}(t) = w(t)$ a.e. $t \in [0,T]$ and $x^*(\cdot)$ satisfies the transversality conditions (*c*),(*d*).

**Theorem 3.1** *Let* $F(\cdot, t): \mathbb{R}^n \rightrightarrows \mathbb{R}^n$ *be a convex mapping and* $f: \mathbb{R}^{2n} \to \mathbb{R}^1 \cup \{+\infty\}$ *be continuous proper convex function. Besides, let* $S \subseteq \mathbb{R}^{2n}$ *be a convex set and* $X: [0,T] \rightrightarrows \mathbb{R}^n$ *be a convex-valued mapping. Suppose that there exists a pair of functions* $\{x^*(\cdot), \tilde{x}(\cdot)\}$ *satisfying a.e. the Euler-Lagrange type inclusions (a), (b) and transversality conditions (c),(d). Then the trajectory* $\tilde{x}(\cdot)$ *is optimal in the convex problem (PHC).*

**Proof.** The Euler-Lagrange inclusion, which is equivalent to subdifferential inclusion

$$(-1)^{(k)} x^{*(k)}(t) \in \partial_x H_F(\tilde{x}(t), x^*(t)) + \upsilon^*(t),\quad \upsilon^*(t) \in K^*_{X(t)}(\tilde{x}(t))$$

implies that

$$H_F(x(t), x^*(t)) - H_F(\tilde{x}(t), x^*(t)) \le \left\langle (-1)^k x^{*(k)}(t) - \upsilon^*(t), x(t) - \tilde{x}(t)\right\rangle. \tag{5}$$

Then from the inequality (5) we have

$$\left\langle x^{(k)}(t)-\tilde{x}^{(k)}(t), x^{*}(t)\right\rangle-\left\langle(-1)^{k} x^{*(k)}(t), x(t)-\tilde{x}(t)\right\rangle \leq 0 . \quad (6)$$

Now we need to integrate inequality (6) over the time interval $[0,T]$:

$$\int_0^T\left[\left\langle x^{(k)}(t)-\tilde{x}^{(k)}(t), x^{*}(t)\right\rangle-\left\langle(-1)^{k} x^{*(k)}(t), x(t)-\tilde{x}(t)\right\rangle\right] dt \leq 0 . \quad (7)$$

Then the square brackets of inequality (7) can be reduced into the following equality relation

$$\begin{aligned}
&\left\langle x^{(k)}(t)-\tilde{x}^{(k)}(t), x^{*}(t)\right\rangle-\left\langle(-1)^{k} x^{*(k)}(t), x(t)-\tilde{x}(t)\right\rangle=-\frac{d}{dt}\left\langle(-1)^{k} x^{*(k-1)}(t), x(t)-\tilde{x}(t)\right\rangle \\
&+\frac{d}{dt}\left\langle(-1)^{k-1} x^{*(k-2)}(t), x'(t)-\tilde{x}'(t)\right\rangle-\frac{d}{dt}\left\langle(-1)^{k-2} x^{*(k-3)}(t), x''(t)-\tilde{x}''(t)\right\rangle \\
&+\frac{d}{dt}\left\langle(-1)^{k-3} x^{*(k-4)}(t), x'''(t)-\tilde{x}'''(t)\right\rangle-\cdots+\frac{d}{dt}\left\langle x^{*}(t), x^{(k-1)}(t)-\tilde{x}^{(k-1)}(t)\right\rangle . \quad (8)
\end{aligned}$$

Now, if we integrate inequality (8) over $[0,T]$ according to higher-order differential calculus [19], we obtain

$$\begin{aligned}
&\int_0^T\left[\left\langle x^{(k)}(t)-\tilde{x}^{(k)}(t), x^{*}(t)\right\rangle-\left\langle(-1)^{k} x^{*(k)}(t), x(t)-\tilde{x}(t)\right\rangle\right] dt \\
&=-\left\langle x(T)-\tilde{x}(T),(-1)^{k} x^{*(k-1)}(T)\right\rangle-\left\langle x'(T)-\tilde{x}'(T),(-1)^{k-1} x^{*(k-2)}(T)\right\rangle \\
&-\left\langle x''(T)-\tilde{x}''(T),(-1)^{k-2} x^{*(k-3)}(T)\right\rangle-\cdots \\
&+\left\langle x^{(k-1)}(T)-\tilde{x}^{(k-1)}(T), x^{*}(T)\right\rangle+\left\langle x(0)-\tilde{x}(0),(-1)^{k} x^{*(k-1)}(0)\right\rangle \\
&+\left\langle x'(0)-\tilde{x}'(0),(-1)^{k-1} x^{*(k-2)}(0)\right\rangle+\left\langle x''(0)-\tilde{x}''(0),(-1)^{k-2} x^{*(k-3)}(0)\right\rangle \\
&+\cdots-\left\langle x^{(k-1)}(0)-\tilde{x}^{(k-1)}(0), x^{*}(0)\right\rangle \\
&=\sum_{j=0}^{k-1}\left\langle(-1)^{j+1} x^{*(j)}(0),\ x^{(k-j-1)}(0)-\tilde{x}^{(k-j-1)}(0)\right\rangle \quad (9) \\
&-\sum_{j=0}^{k-1}\left\langle(-1)^{j+1} x^{*(j)}(T),\ x^{(k-j-1)}(T)-\tilde{x}^{(k-j-1)}(T)\right\rangle .
\end{aligned}$$

Denoting right hand side of (9) by $\Omega$ for all feasible $x(\cdot),\ \tilde{x}(\cdot)$, we have

$$\begin{aligned}
\Omega=&\sum_{j=0}^{k-1}\left\langle(-1)^{j+1} x^{*(j)}(0),\ x^{(k-j-1)}(0)-\tilde{x}^{(k-j-1)}(0)\right\rangle \\
&-\sum_{j=0}^{k-1}\left\langle(-1)^{j+1} x^{*(j)}(T),\ x^{(k-j-1)}(T)-\tilde{x}^{(k-j-1)}(T)\right\rangle
\end{aligned}$$

$$= \sum_{j=0}^{k-2} \left\langle (-1)^{j+1} x^{*(j)}(0),\ x^{(k-j-1)}(0) - \tilde{x}^{(k-j-1)}(0) \right\rangle + \left\langle (-1)^k x^{*(k-1)}(0),\ x(0) - \tilde{x}(0) \right\rangle$$

$$-\sum_{j=0}^{k-2} \left\langle (-1)^{j+1} x^{*(j)}(T),\ x^{(k-j-1)}(T) - \tilde{x}^{(k-j-1)}(T) \right\rangle - \left\langle (-1)^k x^{*(k-1)}(T),\ x(T) - \tilde{x}(T) \right\rangle \le 0 \qquad (10)$$

By the transversality condition ($d$) we have

$$\left\langle (-1)^{j+1} x^{*(j)}(0), x^{(k-j-1)}(0) - \tilde{x}^{(k-j-1)}(0) \right\rangle$$

$$+ \left\langle (-1)^j x^{*(j)}(T), x^{(k-j-1)}(T) - \tilde{x}^{(k-j-1)}(T) \right\rangle \ge 0,\ j = 0,...,k-2. \qquad (11)$$

Then from (10) and (11) for all feasible trajectories $x(\cdot),\ \tilde{x}(\cdot)$ as a result we derive

$$\left\langle (-1)^k x^{*(k-1)}(0),\ x(0) - \tilde{x}(0) \right\rangle - \left\langle (-1)^k x^{*(k-1)}(T),\ x(T) - \tilde{x}(T) \right\rangle \le 0. \qquad (12)$$

On the other hand, let a pair $\left(\mu^*(0), \mu^*(T)\right)$ be an element of the dual cone in the transversality condition ($c$), i.e., $\left(\mu^*(0), \mu^*(T)\right) \in K_S^*(\tilde{x}(0), \tilde{x}(T))$. Then by the condition ($c$) we have

$$f(x(0), x(T)) - f(\tilde{x}(0), \tilde{x}(T)) \ge \left\langle (-1)^{k-1} x^{*(k-1)}(0) + \mu^*(0), x(0) - \tilde{x}(0) \right\rangle$$

$$+ \left\langle \mu^*(T) + (-1)^k x^{*(k-1)}(T), x(T) - \tilde{x}(T) \right\rangle$$

whence

$$f(x(0), x(T)) - f(\tilde{x}(0), \tilde{x}(T)) \ge \left\langle (-1)^{k-1} x^{*(k-1)}(0), x(0) - \tilde{x}(0) \right\rangle$$

$$+ \left\langle (-1)^k x^{*(k-1)}(T), x(T) - \tilde{x}(T) \right\rangle \qquad (13)$$

Then from (12) and (13) for all for all feasible $x(\cdot)$ we have the needed inequality

$$f(x(0), x(T)) - f(\tilde{x}(0), \tilde{x}(T)) \ge 0, \text{ i. e., } f(x(0), x(T)) \ge f(\tilde{x}(0), \tilde{x}(T)). \qquad \square$$

**Remark 3.1** Note that if $D^k$ is an operator of derivatives of the $k$-th order, then the operator $A$, defined as $Ax = D^k x$, is either a self-adjoint operator or an anti-self-adjoint operator depending on the parity of $k$, i.e., if $k$ is an even number, then using the formal adjoint definition of the adjoint operator, we have $A^* x^* := (-1)^k D^k x^* = D^k x^* = Ax^*$, and if $k$ is odd, then $A^* x^* := (-1)^k D^k x^* = -D^k x^* = -Ax^*$. In fact, as can be seen from the proof of the Theorem 3.1, it is also valid in the case $k = 1$ for a problem (PHC) with a first-order differential inclusion. Obviously, the Euler - Lagrange inclusion with a first-order anti-self-adjoint operator and the transversality condition have the form

(1) $-x^{*\prime}(t) \in F^*\left(x^*(t); (\tilde{x}(t), \tilde{x}'(t)), t\right) + K^*_{X(t)}(\tilde{x}(t)),$, a.e. $t \in [0, T]$,

(2) $\tilde{x}'(t) \in F_A\left(\tilde{x}(t); x^*(t), t\right)$,

(3) $\left(x^*(0), -x^*(T)\right) \in \partial f(\tilde{x}(0), \tilde{x}(T)) - K_S^*(\tilde{x}(0), \tilde{x}(T))$.

In this sense the Euler-Lagrange inclusion of problem (PHC) is an immediate generalization of classical Euler-Lagrange inclusion for first order DFIs. $\square$

At the end of paper consider problem (PHC) for third order "linear" DFIs. Hence it follows that inclusion (a) of Theorem 3.1 is a generalization of the Euler - Lagrange inclusion.
Let us consider the problem:

$$\text{minimize } f(x(0),x(T)),$$

(PTL) $$\frac{d^3x(t)}{dt^3} \in F(x(t)), \text{ a.e. } t\in[0,T], \; F(x)\equiv Ax+BU,$$

$$\left(x^{(j)}(0),x^{(j)}(T)\right)\in S,\; j=0,1,2;\quad x(t)\in X(t),\;\forall t\in[0,T], \tag{14}$$

where $f$ is continuously differentiable function, $A$ is $n\times n$ matrix and $B$ is $n\times r$ matrix, $U$ -convex compact in $\mathbb{R}^r$. The problem is to find a control function $\tilde{u}(t)\in U$ so that the corresponding solution $\tilde{x}(t)$ minimizes $f(x(0),x(T))$.

**Theorem 3.2** *The arc $\tilde{x}(t)$ according to the control function $\tilde{u}(t)$ is a solution of the problem (PTL), if there exists an absolutely continuous function $x^*(t)$, satisfying the third order adjoint equation, the transversality condition, and the Pontryagin maximum principle:*

$$-x^{*'''}(t)\in A^*x^*(t)+K^*_{X(t)}(\tilde{x}(t)),$$

$$\left(x^{*''}(0),-x^{*''}(T)\right)\in f'(\tilde{x}(0),\tilde{x}(T))-K_S^*(\tilde{x}(0),\tilde{x}(T)),$$

$$\left((-1)^{j+1}x^{*(j)}(0),(-1)^j x^{*(j)}(T)\right)\in K_S^*\left(\tilde{x}^{(2-j)}(0),\tilde{x}^{(2-j)}(T)\right),\; j=0,1$$

$$\left\langle B\tilde{u}(t),\, x^*(t)\right\rangle=\max_{u\in U}\left\langle Bu,\, x^*(t)\right\rangle.$$

**Proof.** By elementary computations, we find that if $\tilde{v}=A\tilde{x}+B\tilde{u}$, then

$$F^*(v^*;(\tilde{x},\tilde{v}))=\begin{cases} A^*v^*, & \text{if } -B^*v^*\in K_U^*(\tilde{u}),\\ \varnothing, & \text{if } -B^*v^*\notin K_U^*(\tilde{u}),\end{cases}$$

whereas $\left\langle u-\tilde{u},-B^*v^*\right\rangle\ge 0,\, u\in U$ or $\left\langle B\tilde{u},v^*\right\rangle=\max_{u\in U}\left\langle Bu,v^*\right\rangle$. Thus, using Theorem 3.1 we deduce the adjoint linear differential equation of the third order, the transversality conditions, and the Pontryagin maximum principle:

$$-x^{*'''}(t)\in A^*x^*(t)+K^*_{X(t)}(\tilde{x}(t)),\,,$$

$$\left(x^{*''}(0),-x^{*''}(T)\right)\in f'(\tilde{x}(0),\tilde{x}(T))-K_S^*(\tilde{x}(0),\tilde{x}(T)),$$

$$\left((-1)^{j+1}x^{*(j)}(0),(-1)^j x^{*(j)}(T)\right)\in K_S^*\left(\tilde{x}^{(2-j)}(0),\tilde{x}^{(2-j)}(T)\right),\; j=0,1,$$

$$\left\langle B\tilde{u}(t),\, x^*(t)\right\rangle=\max_{u\in U}\left\langle Bu,\, x^*(t)\right\rangle.$$

The proof is completed. □

## 4. The duality to *k*-th order DFIs

We call the following problem, labelled (PHC*), the dual problem to the primary continuous convex problem (PHC):

(PHC*) $$\sup J^*\left[\upsilon^*(\cdot),x^{*(k)}(\cdot);\mu^*(t);x^{*(j)}(t),t=0,T,\,j=0,\dots,k-1\right]$$

where $J*$ is defined as follows:

$$J*\left[\upsilon^*(\cdot),x^{*(k)}(\cdot);\mu^*(t),x^{*(j)}(t),t=0,T,\ j=0,...,k-1\right]=$$

$$-f*\left((-1)^{k-1}x^{*(k-1)}(0)+\mu^*(0),\mu^*(T)+(-1)^{k-1}x^{*(k-1)}(T)\right)$$

$$+\int_0^T M_F\left((-1)^k x^{*(k)}(t)-\upsilon^*(t),x^*(t)\right)dt-\int_0^T W_{X(t)}(-\upsilon^*(t))dt-W_S(-\mu^*(0),-\mu^*(T))$$

$$-\sum_{j=0}^{k-2}W_S\left((-1)^j x^{*(j)}(0),(-1)^{j+1}x^{*(j)}(T)\right).$$

Note that here maximization is carried out over a set of functions $\Re_{k,j}(T,0)\equiv\left\{\upsilon^*(\cdot),x^{*(k)}(\cdot);\mu^*(t),\ x^{*(j)}(t),\ t=0,T,\ j=0,...,k-1\right\}$. Then, according to this notation $\tilde{\Re}_{k,j}(T,0)\equiv\left\{\tilde{\upsilon}^*(\cdot),\tilde{x}^{*(k)}(\cdot);\tilde{\mu}^*(t),\ \tilde{x}^{*(j)}(t),t=0,T,\ j=0,...,k-1\right\}$. We suppose that $\upsilon^*(t),\ t\in[0,T]$ is absolutely continuous function and $x^*(\cdot)\in AC^{k-1}\left(\left[0,T\right],\mathbb{R}^n\right),\ x^{*(k)}(\cdot)\in L^1\left(\left[0,T\right],\mathbb{R}^n\right)$.

Below we present our duality theorem for problem (PHC).

**Theorem 4.1** *Let* $\tilde{x}(t)$ *be the optimal solution to the convex problem (PHC). Then the collection* $\tilde{\Re}_{k,j}(T,0)$ *is a solution to the dual problem (PHC*) if and only if conditions (a) -(d) of Theorem 3.1 are satisfied. In addition, the optimal values in problems (PHC) and (PHC*) are equal.*

**Proof.** First, we show that for all feasible solutions $x(\cdot)$ and a collection of dual variables $\Re_{k,j}(T,0)\equiv\left\{\upsilon^*(\cdot),x^{*(k)}(\cdot);\mu^*(t),\ t=0,T,\ j=0,...,k-1\right\}$, the value of the problem (PHC) is not less than the value (PCH *), that is

$$f(x(0),x(T))\geq -f*\left((-1)^{k-1}x^{*(k-1)}(0)+\mu^*(0),\mu^*(T)+(-1)^{k-1}x^{*(k-1)}(T)\right)$$

$$+\int_0^T M_F\left((-1)^k x^{*(k)}(t)-\upsilon^*(t),x^*(t)\right)dt-\int_0^T W_{X(t)}(-\upsilon^*(t))dt-W_S(-\mu^*(0),-\mu^*(T))$$

$$-\sum_{j=0}^{k-2}W_S\left((-1)^j x^{*(j)}(0),(-1)^{j+1}x^{*(j)}(T)\right). \qquad (15)$$

Clearly, applying the Young's inequality [18], we can write

$$-f*\left((-1)^{k-1}x^{*(k-1)}(0)+\mu^*(0),\mu^*(T)+(-1)^{k-1}x^{*(k-1)}(T)\right)\leq f(x(0),x(T))$$

$$-\left\langle x(0),(-1)^{k-1}x^{*(k-1)}(0)+\mu^*(0)\right\rangle-\left\langle x(T),\mu^*(T)+(-1)^{k-1}x^{*(k-1)}(T)\right\rangle. \qquad (16)$$

Now, using the definitions of $M_F$ and Hamiltonian functions, we have

$$M_F\left((-1)^k x^{*(k)}(t)-\upsilon^*(t),\ x^*(t)\right)\leq\left\langle x(t),(-1)^k x^{*(k)}(t)-\upsilon^*(t)\right\rangle \qquad (17)$$

$$-\left\langle x^{(\kappa)}(t),x^*(t)\right\rangle=\left\langle x(t),(-1)^k x^{*(k)}(t)\right\rangle-\left\langle x^{(\kappa)}(t),x^*(t)\right\rangle-\left\langle x(t),\upsilon^*(t)\right\rangle.$$

Hence, integrating (17) over interval $\left[0,T\right]$ we have

$$\int_0^T M_F\left((-1)^k x^{*(k)}(t) - \upsilon^*(t),\ x^*(t)\right)dt$$

$$\le \int_0^T \left[\left\langle x(t), (-1)^\kappa x^{*(\kappa)}(t)\right\rangle - \left\langle x^{(\kappa)}(t), x^*(t)\right\rangle\right]dt - \int_0^T \left\langle x(t), \upsilon^*(t)\right\rangle dt\,. \tag{18}$$

On the other hand, it can be easily seen that

$$-\int_0^T W_{X(t)}(-\upsilon^*(t))dt \le \int_0^T \left\langle x(t), \upsilon^*(t)\right\rangle dt \tag{19}$$

$$-W_S(-\mu^*(0), -\mu^*(T)) \le \left\langle x(0), \mu^*(0)\right\rangle + \left\langle x(T), \mu^*(T)\right\rangle, \tag{20}$$

$$-W_S\left((-1)^j x^{*(j)}(0), (-1)^{j+1} x^{*(j)}(T)\right)$$

$$\le \left\langle (-1)^{j+1} x^{*(j)}(0),\ x^{(k-j-1)}(0)\right\rangle + \left\langle (-1)^j x^{*(j)}(T),\ x^{(k-j-1)}(T)\right\rangle,\ j = 0,\dots,k-2\,. \tag{21}$$

Summing now the inequalities (16) and (18)-(21) we derive

$$J*\left[\upsilon^*(\cdot), x^{*(k)}(\cdot), \mu^*(t), x^{*(j)}(t), t = 0, T;\ j = 0,\dots,k-1\right]$$

$$\le f(x(0), x(T)) - \left\langle x(0), (-1)^{k-1} x^{*(k-1)}(0)\right\rangle$$

$$-\left\langle x(T), (-1)^{k-1} x^{*(k-1)}(T)\right\rangle + \int_0^T \left[\left\langle x(t), (-1)^\kappa x^{*(\kappa)}(t)\right\rangle - \left\langle x^{(\kappa)}(t), x^*(t)\right\rangle\right]dt$$

$$+\sum_{j=0}^{k-2}\left[\left\langle (-1)^{j+1} x^{*(j)}(0),\ x^{(k-j-1)}(0)\right\rangle + \left\langle (-1)^j x^{*(j)}(T),\ x^{(k-j-1)}(T)\right\rangle\right]. \tag{22}$$

We need calculate the integral of dot product differences $Q = \left\langle x(t), (-1)^\kappa x^{*(\kappa)}(t)\right\rangle$ $-\left\langle x^{(\kappa)}(t), x^*(t)\right\rangle$. To this end we transform $Q$ as follows

$$Q = \frac{d}{dt}\left\langle x(t), (-1)^\kappa x^{*(\kappa-1)}(t)\right\rangle + \frac{d}{dt}\left\langle x'(t), (-1)^{\kappa-1} x^{*(\kappa-2)}(t)\right\rangle$$

$$+\frac{d}{dt}\left\langle x''(t), (-1)^{\kappa-2} x^{*(\kappa-3)}(t)\right\rangle + \cdots + \frac{d}{dt}\left\langle x^{(\kappa-2)}(t), x^{*\prime}(t)\right\rangle - \frac{d}{dt}\left\langle x^{(\kappa-1)}(t), x^*(t)\right\rangle. \tag{23}$$

Hence, using (23) we can compute the integral of $Q$ over time interval $[0, T]$ as follows:

$$\int_0^T Q dt = \left\langle x(T), (-1)^\kappa x^{*(\kappa-1)}(T)\right\rangle + \left\langle x'(T), (-1)^{\kappa-1} x^{*(\kappa-2)}(T)\right\rangle$$

$$+\left\langle x''(T), (-1)^{\kappa-2} x^{*(\kappa-3)}(T)\right\rangle + \cdots + \left\langle x^{(\kappa-2)}(T), x^{*\prime}(T)\right\rangle$$

$$-\left\langle x^{(\kappa-1)}(T), x^*(T)\right\rangle - \left\langle x(0), (-1)^\kappa x^{*(\kappa-1)}(0)\right\rangle - \left\langle x'(0), (-1)^{k-1} x^{*(k-2)}(0)\right\rangle$$

$$-\left\langle x''(0),(-1)^{k-2}x^{*(k-3)}(0)\right\rangle-\cdots-\left\langle x^{(k-2)}(0),x^{*\prime}(0)\right\rangle+\left\langle x^{(k-1)}(0),x^{*}(0)\right\rangle$$

$$=-\sum_{j=0}^{k-1}\left\langle x^{(k-j-1)}(T),(-1)^{j}x^{*(j)}(T)\right\rangle-\sum_{j=0}^{k-1}\left\langle x^{(k-j-1)}(0),(-1)^{j+1}x^{*(j)}(0)\right\rangle$$

$$=-\sum_{j=0}^{k-1}\left[\left\langle x^{(k-j-1)}(0),(-1)^{j+1}x^{*(j)}(0)\right\rangle+\left\langle x^{(k-j-1)}(T),(-1)^{j}x^{*(j)}(T)\right\rangle\right]. \quad (24)$$

Taking into account (24) in (22) we have

$$J^*\left[\upsilon^*(\cdot),x^{*(k)}(\cdot);\mu^*(t),x^{*(j)}(t),t=0,T,\ j=0,...,k-1\right]$$

$$\le f(x(0),x(T))-\left\langle x(0),(-1)^{k-1}x^{*(k-1)}(0)\right\rangle-\left\langle x(T),(-1)^{k-1}x^{*(k-1)}(T)\right\rangle$$

$$+\sum_{j=0}^{k-2}\left[\left\langle(-1)^{j+1}x^{*(j)}(0),\ x^{(k-j-1)}(0)\right\rangle+\left\langle(-1)^{j}x^{*(j)}(T),\ x^{(k-j-1)}(T)\right\rangle\right]$$

$$-\sum_{j=0}^{k-1}\left[\left\langle x^{(k-j-1)}(0),(-1)^{j+1}x^{*(j)}(0)\right\rangle+\left\langle x^{(k-j-1)}(T),(-1)^{j}x^{*(j)}(T)\right\rangle\right]$$

$$=-\left\langle x(T),(-1)^{k}x^{*(k-1)}(T)\right\rangle-\left\langle x(0),(-1)^{k}x^{*(k-1)}(0)\right\rangle=f(x(0),x(T)). \quad (25)$$

Thus, it follows from (25) that for an arbitrary feasible solution $x(\cdot)$ of problem (PHC) and a family of dual variables $\Re_{k,j}(T,0)$ the inequality holds:

$$J^*\left[\upsilon^*(\cdot),x^{*(k)}(\cdot);\mu^*(t),x^{*(j)}(t),t=0,T,\ j=0,...,k-1\right]\le f(x(0),x(T)).$$

Thus, the justification of inequality (15) is proved.

In addition, assume that the set of dual variables $\tilde{\Re}_{k,j}(T,0)$ satisfies conditions (*a*) - (*d*) of Theorem 3.1. Then, as we saw in inequality (5), the proof of Theorem 3.1, it follows from the Euler-Lagrange type inclusion (*a*) and condition (*b*) that

$$H_F(x(t),\tilde{x}^*(t))-H_F(\tilde{x}(t),\tilde{x}^*(t))\le\left\langle(-1)^k\tilde{x}^{*(k)}(t)-\upsilon^*(t),x(t)-\tilde{x}(t)\right\rangle.$$

In turn, recall that by Lemma 2.6 [18, p.64] (see, Section 2) in order for $x^*$ to be an element of LAM $F^*$, it is necessary and sufficient that $M_F(x^*,v^*)=\left\langle x,x^*\right\rangle-H_F(x,v^*)$, whence

$$\left\langle(-1)^k\tilde{x}^{*(k)}(t)-\tilde{\upsilon}^*(t),\tilde{x}(t)\right\rangle-H_F\left(\tilde{x}(t),\tilde{x}^*(t)\right)=M_F\left((-1)^k\tilde{x}^{*(k)}(t)-\tilde{\upsilon}^*(t),\ \tilde{x}^*(t)\right). \quad (26)$$

Besides, by Theorem 1.27 [18] the transversality conditions (*c*), (*d*) of Theorem 3.1 inscribed for the family $\tilde{\Re}_{k,j}(T,0)$ mean that the inequalities (17)-(21) are satisfied as the exact equalities. Then, the inequality (15) is fulfilled as the equality and for $\tilde{x}(\cdot)$ and $\tilde{\Re}_{k,j}(T,0)$ the equality of values of (PHC) and (PHC*) problems is guaranteed. Thus, the conditions (*a*)-(*d*) of Theorem 3.1 implies that $\tilde{\Re}_{k,j}(T,0)$ is a solution of the dual problem (PHC*). The converse is proved similarly. Since by Lemma 2.6 [18, p.64] $M_F(x^*,v^*)=\left\langle x,x^*\right\rangle-H_F(x,v^*)$ the last formula (26) for our problem means that (17) is satisfied, whence we immediately have an inclusion of Euler-Lagrange type (*a*) of Theorem 3.1. Besides, since the LAM $F^*$ is nonempty, the condition (*b*) of Theorem 3.1 is satisfied. Note that, by assumption, $\tilde{\Re}_{k,j}(T,0)$ is a solution to the dual problem and therefore (16) is fulfilled as an equality, without the transversality conditions. The proof of theorem is completed. □

## 5. Duality problems for third order linear and fourth order polyhedral DFIs

### 5.1 Linear problem

In this subsection, we will construct the problem dual to the continuous problem (PTL) from Section 3. First, we compute $M_F$ function

$$M_F(x^*,v^*)=\inf_{(x,v)\in \mathrm{gph}F}\left\{\left\langle x,x^*\right\rangle-\left\langle v,v^*\right\rangle\right\}$$

$$=\inf_x\left[\left\langle x,x^*-A^*v^*\right\rangle\right]-\sup_{u\in U}\left\langle u,B^*v^*\right\rangle=\begin{cases}-W_U(B^*v^*), & \text{if } x^*=A^*v^*,\\ -\infty, & \text{otherwise.}\end{cases} \tag{27}$$

Then according to the dual problem (PHC*) from (27), we derive

$$M_F\left(-x^{*'''}(t)-\upsilon^*(t),\ x^*(t)\right)=\begin{cases}-W_U(B^*x^*(t)), & \text{if } -x^{*'''}(t)-\upsilon^*(t)=A^*x^*(t),\\ -\infty, & \text{otherwise,}\end{cases}$$

whence we deduce the Euler-Lagrange type adjoint inclusion (equation)

$$\frac{d^3x^*(t)}{dt^3}=-A^*x^*(t)-\upsilon^*(t),\quad \left\langle B\tilde{u}(t),\ x^*(t)\right\rangle=\max_{u\in U}\left\langle Bu,\ x^*(t)\right\rangle. \tag{28}$$

Thus, it can be easily seen that the dual problem to problem (PTL) is

$$\text{(PTL*)}\qquad \sup\Big\{-f^*\left(x^{*''}(0)+\mu^*(0),\mu^*(T)+x^{*''}(T)\right)-\int_0^T W_U(B^*x^*(t))dt$$
$$-\int_0^T W_{X(t)}\left(-x^{*'''}(t)-A^*x^*(t)\right)dt-W_S(-\mu^*(0),-\mu^*(T))$$
$$-\sum_{j=0}^{1}W_S\left((-1)^j x^{*(j)}(0),(-1)^{j+1}x^{*(j)}(T)\right)\Big\}.$$

Therefore, it is interesting to note that maximization in the problem (PTL*) is performed on the set of functions $\Re_{k,j}(T,0)\equiv\left\{x^{*'''}(\cdot);\ \mu^*(t),\ x^{*(j)}(t),\ t=0,T,\ j=0,1,2\right\}$.

**Theorem 5.1** *Let* $\tilde{x}(t)$ *be the optimal solution to the convex problem (PTL). Then the collection* $\tilde{\Re}_{k,j}(T,0)$ *is a solution to the dual problem (PTL*) if and only if conditions (a) -(d) of Theorem 3.2 are satisfied. In addition, the optimal values in problems (PHC) and (PHC *) are equal.*

### 5.2 Polyhedral problem

Here we construct a dual problem (PFC*) to a problem with a fourth-order polyhedral differential inclusion with the endpoint and state conditions

$$\text{infimum } f(x(0),x(T)),$$

$$\text{(PFC)}\qquad \frac{d^4x(t)}{dt^4}\in F\left(x(t)\right),\ \text{a.e. } t\in[0,T],\ F(x)=\left\{v:Ax-Ev\le d\right\},$$

$$\left(x^{(j)}(0),x^{(j)}(T)\right)\in S,\ j=0,1,2,3;,\quad x(t)\in X(t),\ \forall t\in[0,T],$$

where $A, E$ are $m \times n$ dimensional matrices, $d$ is a $m$ -dimensional column-vector, $f(\cdot,\cdot)$ is a proper convex function. The problem is to find the trajectory $\tilde{x}(\cdot)$ of the problem (PFC) that minimizes the Mayer functional $f(\cdot,\cdot)$.

Thus, based on Theorem 3.1 for the problem (PFC), we prove the following theorem.

**Theorem 5.1.** *For the optimality of the trajectory $\tilde{x}(\cdot)$ in problem (PFC) with a fourth-order polyhedral differential inclusion and endpoint and state conditions it is sufficient that there exists a nonnegative function $\lambda(t) \geq 0,\ t \in [0,T]$ satisfying (1), (2):*

(1) $-x^{*(iv)}(t) = A^*\lambda(t),\ \langle A\tilde{x}(t) - E\tilde{x}^{(iv)}(t) - d, \lambda(t)\rangle = 0,$ a.e. $t \in [0,T]$,

(2) $\left(-x^{*'''}(0), x^{*'''}(T)\right) \in \partial f(\tilde{x}(0), \tilde{x}(T)) - K_S^*(\tilde{x}(0), \tilde{x}(T))$,

$\left((-1)^{j+1}E^*\lambda^{(j)}(0), (-1)^j E^*\lambda^{(j)}(T)\right) \in K_S^*\left(\tilde{x}^{(3-j)}(0), \tilde{x}^{(3-j)}(T)\right),\ j = 0,1,2$.

**Proof.** By the condition (1) of Theorem 3.1 one has

$$-x^{*(iv)}(t) \in F^*\left(x^*(t); (\tilde{x}(t), \tilde{x}^{(iv)}(t))\right). \tag{29}$$

Hence, we need calculate the LAM $F^*\left(x^*(t); (\tilde{x}(t), \tilde{x}'''(t))\right)$. Clearly, in this case gph$F$ $= \{(x,v) : Ax - Ev \leq d\}$. For a point $\tilde{w} = (\tilde{x}, \tilde{v}) \in$ gph$F$ we put $I(\tilde{w}) = \{\, i : A_i\tilde{x} - E_i\tilde{v} = d_i,$ $i = 1,\dots,m\}$, where $A_i, E_i$ be the *i*-th row of the matrices $A, E$ respectively, and $d_i$ be the *i*-th component of the vector $d$. On the definition of cone of tangent directions $K_F(\tilde{w}) = \{\bar{w}:$ $\tilde{w} + \gamma\bar{w} \in$ gph$F$ for sufficiently small $\gamma > 0\}$. For each $i \in I(\tilde{w})$ the inequality $A_i(\tilde{x} + \gamma\bar{x})$ $-E_i(\tilde{v} + \gamma\bar{v}) = d_i + \gamma(A_i\bar{x} - E_i\bar{v}) \leq d_i$ is satisfied, if $A_i\bar{x} - E_i\bar{v} \leq 0,\ i \in I(\tilde{w})$. If $i \notin I(\tilde{w})$, then the inequality $A_i(\tilde{x} + \gamma\bar{x}) - E_i(\tilde{v} + \gamma\bar{v}) = (A_i\tilde{x} - E_i\tilde{v}) + \gamma(A_i\bar{x} - E_i\bar{v}) < d_i$ is valid for sufficiently small $\gamma > 0$ regardless of the choice of $\bar{w} = (\bar{x}, \bar{v})$. It follows that $K_F(\tilde{w}) = \{\bar{w} : A_i\bar{x} - E_i\bar{v}) \leq 0,$ $i \in I(\tilde{w}\}$. According to Farkas theorem [18, p. 22], it is not hard to see that $(x^*, v^*) \in K_F^*(\tilde{w})$ if and only if

$$x^* = -\sum_{i \in I(\tilde{w})} A_i^*\lambda_i,\ v^* = \sum_{i \in I(\tilde{w})} E_i^*\lambda_i, \lambda_i \geq 0, i = 1,\dots,m, \tag{30}$$

where $A_i^*,\ E_i^*$ are vector-columns. Then, setting $\lambda_i = 0, i \notin I(\tilde{w})$ and denoting by $\lambda$ the vector-column with components $\lambda_i$ we derive from (30) that

$$K_F^*(\tilde{w}) = \left\{(x^*, v^*) : x^* = -A^*\lambda,\ v^* = E^*\lambda,\ \lambda \geq 0,\ \langle A\tilde{x} - E\tilde{v} - d, \lambda\rangle = 0\right\}.$$

Finally, for the polyhedral LAM from (30) we have the following formula

$$F^*\left(v^*; (\tilde{x}, \tilde{v})\right) = \{-A^*\lambda :\ v^* = E^*\lambda,\ \lambda \geq 0,\ \langle Ax - Ev - d, \lambda\rangle = 0\}. \tag{31}$$

In fact, from (30),(31) it should be noted that $F^*\left(v^*; (\tilde{x}, \tilde{v})\right)$ does not depend on point $\tilde{w} = (\tilde{x}, \tilde{v})$, but depends on the set $I(\tilde{w})$ (since the number of such sets is finite it follows that the number of different LAM $F^*(v^*; (\tilde{x}, \tilde{v}))$ is finite). Thus, from (29) and (31) we derive that

$$-x^{*(iv)}(t) = A^*\lambda(t),\ \langle A\tilde{x}(t) - E\tilde{x}^{(iv)}(t) - d, \lambda(t)\rangle = 0,\ \text{a.e. } t \in [0,T], \tag{32}$$

where $x^*(t) = E^*\lambda(t)$.

Therefore, since $x^*(t)=E*\lambda(t)$ the transversality condition of Theorem 3.1

$$\left(-E*\lambda'''(0),E*\lambda'''(T)\right)\in\partial f(\tilde{x}(0),\tilde{x}(T))-K_S*(\tilde{x}(0),\tilde{x}(T)),$$

$$\left((-1)^{j+1}x^{*(j)}(0),(-1)^j x^{*(j)}(T)\right)\in K_S*\left(\tilde{x}^{(3-j)}(0),\tilde{x}^{(3-j)}(T)\right),\ j=0,1,2,$$

has the following form

$$\left(-x^{*'''}(0),x^{*'''}(T)\right)\in\partial f(\tilde{x}(0),\tilde{x}(T))-K_S*(\tilde{x}(0),\tilde{x}(T)),$$

$$\left((-1)^{j+1}E*\lambda^{(j)}(0),(-1)^j E*\lambda^{(j)}(T)\right)\in K_S*\left(\tilde{x}^{(3-j)}(0),\tilde{x}^{(3-j)}(T)\right),\ j=0,1,2. \qquad \square$$

It remains to construct the dual problem (PFC*) to problem (PFC). First of all, according to the dual problem (PHC*) we should compute the $M_F$ function:

$$M_F(x^*,v^*)=\inf\left\{\langle x,x^*\rangle-\langle v,v^*\rangle:(x,v)\in \text{gph}F\right\}.$$

It can be easily seen that, denoting $w=(x,v)\in\mathbb{R}^{2n}, w^*=(x^*,-v^*)\in\mathbb{R}^{2n}$ we have a problem

$$\inf\left\{\langle w,w^*\rangle:Cw\le d\right\}, \tag{33}$$

where $C=\left[A\,\vdots -E\right]$ is $m\times 2n$ block matrix. Then for a solution $\tilde{w}=(\tilde{x},\tilde{v})$ of (33) there exists $m$-dimensional vector $\lambda\ge 0$ such that $w^*=-C*\lambda$, $\langle A\tilde{x}-E\tilde{v}-d,\lambda\rangle=0$. Thus, $w^*=-C*\lambda$ implies that $x^*=-A*\lambda$, $v^*=-E*\lambda,\lambda\ge 0$. Then, we find that

$$M_F(x^*,v^*)=\langle\tilde{x},-A*\lambda\rangle-\langle\tilde{v},-E*\lambda\rangle=-\langle A\tilde{x},\lambda\rangle+\langle E\tilde{v},\lambda\rangle=-\langle d,\lambda\rangle. \tag{34}$$

Besides, taking into account the function $M_F\left(x^{*(iv)}(t)-\upsilon^*(t),x^*(t)\right)$ and the first and second relations (34), we derive

$$x^{*(iv)}(t)=-A*\lambda(t)+\upsilon^*(t),\ x^*(t)=-E*\lambda(t),\ \lambda(t)\ge 0,$$

whence

$$-E*\lambda^{(iv)}(t)=-A*\lambda(t)+\upsilon^*(t),\ \lambda(t)\ge 0. \tag{35}$$

Consequently, taking into account (34), (35), and duality Theorem 4.1 under the conditions (2) of Theorem 5.1 we have the following dual problem for fourth order Polyhedral DFIs

(PFC*)

$$\sup\Bigg\{-f*\left(-x^{*'''}(0)+\mu^*(0),\mu^*(T)-x^{*'''}(T)\right)-\int_0^T\langle d,\lambda(t)\rangle dt$$
$$-\int_0^T W_{X(t)}\left(E*\lambda^{(iv)}(t)-A*\lambda(t)\right)dt-W_S(-\mu^*(0),-\mu^*(T))$$
$$-\sum_{j=0}^{2}W_S\left((-1)^{j+1}E*\lambda^{(j)}(0),(-1)^j E*\lambda^{(j)}(T)\right)\Bigg\}.$$

Here maximization in this problem is realised on the set of functions $\Re_{k,j}(T,0)\equiv\{\lambda^{(iv)}(\cdot);$ $\mu^*(t),\ \lambda^{(j)}(t),\ t=0,T;\ j=0,1,2,3\}$.

Now, based on Theorems 3.1 for problem (PFC), we have

**Theorem 5.2** *Let $\tilde{x}(t)$ be the optimal solution to the convex problem (PFC). Then the collection $\{\tilde{\lambda}^{(iv)}(\cdot)\ge 0,\tilde{\mu}^*(t),\ \tilde{\lambda}^{(j)}(t),\ t=0,T;\ j=0,1,2,3\}$, $t\in[0,T]$ is a solution to the dual problem (PFC*) if and only if conditions (1),(2) of Theorem 5.1 are satisfied. In addition, the optimal* *values in problems (PFC) and (PFC *) are equal.*